\documentclass[12pt]{article}
\usepackage{amssymb} \usepackage{amsfonts}
\usepackage{amsmath}
\begin{document}
\title{Eta invariants for flat manifolds}
\author{A. Szczepa\'nski
\thanks{The author was supported by Max Planck Institute in Bonn}\\
Institute of Mathematics, University of Gda\'nsk\\
ul. Wita Stwosza 57,
80-952 Gda\'nsk,
Poland\\
E-mail: matas@univ.gda.pl}
\date{\today}
\maketitle
\begin{abstract}
Using a formula from H. Donnelly \cite{HD}, we prove that for a  
family of seven
dimensional flat manifolds with cyclic holonomy groups the $\eta$ invariant
of the signature operator is an integer number. 
We also present an infinite family of flat manifolds
with integral $\eta$ invariant.
Our main motivation
is a paper of D. D. Long and A. Reid  "On the geometric boundaries of hyperbolic $4$-manifolds", see \cite{LR}.  
\end{abstract}
\newcommand{\Z}{{\mathbb Z}}
\newcommand{\Q}{{\mathbb Q}}
\newcommand{\R}{{\mathbb R}}
\newcommand{\C}{{\mathbb C}}
\newcommand{\h}{{\mathbb H}}
\newcommand{\N}{{\mathbb N}}
\section{Introduction}
Let $M^n$ be a closed Riemannian manifold of dimension $n.$
We shall call $M^n$ flat if, at any point, the sectional curvature is equal to zero. Equivalently,$M^n$ is isometric
to the orbit space $\R^n/\Gamma,$ where $\Gamma$ is a discrete, torsion-free and co-compact subgroup of $O(n)\ltimes\R^n$ = Isom($\R^n$). 
From the Bieberbach theorem (see \cite{charlap}, \cite{S}, \cite{wolf}) $\Gamma$ defines the short exact sequence of groups
\begin{equation}\label{crysb}
0\rightarrow \Z^n\rightarrow\Gamma\stackrel{p}\rightarrow G\rightarrow 0,
\end{equation}
where $G$ is a finite group. 
$\Gamma$ is called a Bieberbach group and $G$ its holonomy group.
We can define a holonomy representation $\phi:G\to GL(n,\Z)$ by the formula:
\begin{equation}\label{holonomyrep}
\forall g\in G,\phi(g)(e_i) = \tilde{g}e_i(\tilde{g})^{-1},
\end{equation}
where $e_i\in\Gamma$ are generators of $\Z^n$ for $i=1,2,...,n,$ and $\tilde{g}\in\Gamma$
such that $p(\tilde{g})=g.$

Our main motivation is the paper of D. D. Long and A. W. Reid \cite{LR}. 
Using the methods from \cite{APS}, the authors of \cite{LR} proved that an
obstruction, for the flat $4n-1$-dimensional manifold, to be realized as
the cusp cross-section of a complete 
finite volume one-cusped hyperbolic $4n$-manifold, 
is the non-integrality of the $\eta$ invariant of the signature operator. They gave (see \cite{LR}) an example of a 3-dimensional
flat manifold $M^3$ with $\eta(M^3)\notin \Z,$ see Example 1.

H. Donnelly in \cite{HD} formulated a general formula for the above $\eta$ invariant for some special class of flat manifolds.
From (\ref{crysb}) it is easy to see that any flat manifold $M^n$ is diffeomorphic
to $T^n/G,$ where $T^n = \R^n/\Z^n$ is the $n$-dimensional torus. 
Hence, we can say that a map $$T^n\to T^n/G$$ is regular covering of $T^n/G$ with covering group $G$.
The above formula expresses the $\eta$ invariant of the quotient space $T^n/G$ 
with the $\eta$ invariant
of $T^n,$ and some properties of the covering map (or a group action). 
Such approach was already considered in an original Atiyah, Patodi, Singer paper \cite[pp. 408-413]{APS}. 

Let $T^{4n-2}$ be any flat $(4n-2)$-dimensional torus, $S^1$ be the unit circle and $G$ be a finite
group which acts on $S^1\times T^{4n-2} = T^{4n-1},$ such that 
$T^{4n-1}/G$ is an oriented flat manifold with holonomy group $G.$
Let $\Gamma = \pi_1(T^{4n-1}/G), g\in G$ and
$g = p(\tilde{g}) = \bar{A},$ where $\tilde{g} = (\bar{A},b)\in\Gamma\subset SL(4n-1,\Z)\ltimes\R^{4n-1}.$
We assume that $g$ acts on $T^{4n-1}$ in the following way
\begin{equation}\label{action}
g(x,\bar{x}) = (x+a,A\bar{x}+\bar{a}),
\end{equation}
where $b=(a,\bar{a})\in S^1\times T^{4n-2}$ and $A\in SL(4n-2,\Z).$
Equivalently, it means that
$\bar{A} =
\left[
\begin{smallmatrix}
1&0\\
0&A
\end{smallmatrix}\right].$
Here an element $b\in\R^{4n-1}$
also denotes its image in $\R^{4n-1}/\Z^{4n-1}.$
Since $\Gamma$ is torsion free we can assume that $0\neq a\in \R/\Z$ for $g\neq 1.$
\vskip 1mm
\noindent
Our main result (Theorem 1) is the following:
\vskip 1mm
\noindent
{\em If $M^7$ is a seven dimensional, oriented flat manifold with cyclic holonomy group, 
which satisfies condition (\ref{action}), then $\eta(M^7)\in\Z.$}
\vskip 1mm
\noindent
There exists a classification of flat manifolds up to dimension six, (see \cite{OPS}). 
It was computed with the support of computer system CARAT. 
This algorithm also gives a method for the classification of seven dimensional
Bieberbach groups with a cyclic holonomy group.
It was done by R. Lutowski (see \cite{RL} and \cite{Lu10}).
In the proof of the above result we were concentrated only on those 
oriented flat manifolds which satisfy condition (\ref{action}). 
Summing up, if there exists a seven dimensional 
flat manifold with the $\eta$ invariant $\notin \Z,$ then either it has a 
noncyclic holonomy group
or it has a cyclic holonomy with a special holonomy representation.  
\vskip 1mm
There is already some literature on the $\eta$ invariant of flat manifolds. For example see
\cite{GMP},\cite{MP},\cite{Pf},\cite{SS}. However, in all these articles the authors mainly are
concentrated on the $\eta$ invariant of the Dirac operator. 
\vskip 1mm
\noindent
Let us present a structure of the paper. We prove Theorem 1 in section 3.
For the proof, we use a generalized formula from H. Donnelly \cite{HD}. It is
recalled in section 2, see Propositions \ref{prop1}, \ref{prop2} and Remark 1.
In the last section we present two families of flat manifolds which exactly satisfy 
assumptions from \cite{HD}. 
For the first family
\footnote{So called Hantzsche-Wendt manifolds.}
of oriented $n$-dimensional flat manifolds with the holonomy group $(\Z_2)^{n-1}$
(see \cite{HW}), we prove that the
$\eta$ invariant is always equal to zero. 
In the case of the second family we give an exact formula (\ref{wzor}) for the 
$\eta$ invariant.
However, we do not know how to prove 
that the values of the $\eta$ invariant are in $\Z.$ 
We do it only in a very special case.  
\vskip 3mm
\noindent
For the proof of the main result we use the computer package CARAT, see \cite{OPS}, \cite{RL} and \cite{Lu10}. 
We thank R. Lutowski for his assistance
in the use of CARAT and checking our calculations in the proof of Theorem 1.
Moreover, the author would like to thank W. Miklaszewski, B. Putrycz for their help 
in the use of "MATHEMATICA" version 7 and M. Mroczkowski
for improving English. Finally, we would like thank the referee(s) for a careful reading
and many constructive remarks.
\section{Donnelly's formula}
Since the result of H. Donnelly \cite{HD} is more then 30 years old let us recall it with some details and comments.
We keep the notations from the introduction.
Let $X^{4n} = X$ be a compact oriented Riemannian manifold of dimension $4n$ with non-empty boundary $Y^{4n-1} = Y.$
Assume that the metric of $X$ is a product near the boundary $Y.$
Let $\Lambda(Y)$ be the exterior algebra of $Y$
 (see \cite{APSI}, \cite{HD}), and let
$B:\Lambda(Y)\to\Lambda(Y)$ be a first order self-adjoint elliptic differential operator defined by
$$B\phi = (-1)^{n+p+1}(\epsilon\ast d-d\ast)\phi,$$
where $\ast$ is the duality operator on $Y$ and $\phi$ is either a $2p$-form ($\epsilon =1$)
or a $(2p-1)$-form ($\epsilon = -1$). $B$ preserves the parity of forms on $Y$ and commutes
with $\phi\mapsto (-1)^p\ast\phi,$ so that $B = B^{ev}\oplus B^{odd}$ and $B^{ev}$ is isomorphic to $B^{odd}.$
$B^{ev}$ has pure point spectrum consisting
of eigenvalues $\lambda$ with multiplicity dim($\lambda$). The spectral function
$$\eta(s,Y) = \Sigma_{\lambda\neq 0}(\text{sign}\lambda)(\text{dim}\lambda)\mid\lambda\mid^{-s}$$
converges for Re($s$) sufficiently large and has a meromorphic continuation to the entire complex $s$-plane.
Moreover $\eta(0,Y)$ is finite, see \cite{APSI}.
\vskip 1mm
Consider the finite group $G$ acting isometrically on a manifold $Y$
and suppose $g\in G.$ Then the map defined by $g$ on sections of $\Lambda^{ev}(Y)$ commutes with $B^{ev}.$
This induces linear maps $g_{\lambda}^{\ast}$ on each eigenspace, with eigenvalue $\lambda,$ of $B^{ev}.$
The spectral function
$$\eta_{g}(x,Y) = \Sigma_{\lambda\neq 0}(\text{sign}\lambda)\text{Tr}(g_{\lambda}^{\ast})\mid\lambda\mid^{-s}$$
converges for Re$(s)$ sufficiently large and has a 
meromorphic continuation to the entire complex $s$-plane.
Suppose that $\hat{Y}\to Y$ is a regular covering space with finite covering group $G$ of order $\mid G\mid.$
For each irreducible unitary representation $\alpha$ of $G,$ one has an associated flat vector bundle $E_{\alpha}\to Y.$
The invariants $\eta_{\alpha}(0,Y)$ are defined using the spectrum of the operator
$B^{ev}_{\alpha}:\Lambda^{ev}(Y)\otimes E_{\alpha}\to\Lambda^{ev}(Y)\otimes E_{\alpha}.$ 
These invariants were studied in \cite{APS}. In particular
\begin{equation}\label{general}
\eta_{\alpha}(0,Y) = \tfrac{1}{\mid G\mid}\Sigma_{g\in G}\eta_g(0,\hat{Y})\chi_{\alpha}(g),
\end{equation}
where $\chi_{\alpha}$ is the character of $\alpha.$ We use a special version of formula (\ref{general}).
If we take $\alpha$ to be the trivial one-dimensional representation in (\ref{general}) then we have
\begin{equation}\label{generalI}
\eta(0,\hat{Y}) - \mid G\mid\eta(0,Y) = - \Sigma_{g\neq 1}\eta_{g}(0,\hat{Y}),
\end{equation}
where the sum on the right is taken over group elements $g\in G, g\neq 1.$
Moreover, in our case for $\hat{Y} = T^{4n-1}$ we have $\eta(0,\hat{Y}) = 0,$ cf. \cite[p. 410]{APS}. 
 
\vskip 1mm
\noindent
In \cite{HD} the following is proved:
\newtheorem{prop}{Proposition}
\begin{prop}\label{prop1} {\em (\cite[Proposition 4.6]{HD})}
Let $g:T^{4n-1}\to T^{4n-1}$ be given by the formula {\em (\ref{action})} with $A\in SO(4n-2,\Z)$ and $a\neq 0.$
If $A$ has $1$ as an eigenvalue, then $\eta_{g}(0,T^{4n-1}) = 0.$
\end{prop}  
\vskip 5mm
\newtheorem{rem}{Remark} 
\noindent
Hence Proposition~\ref{prop1} reduces the problem of computing $\eta_{g}(0,T^{4n-1})$ only
to those isometries $g$ satisfying det($I-A$)$\neq 0.$ 
For such $g$ the following is proved in \cite{HD}:
\begin{prop}\label{prop2} {\em (\cite[Proposition 4.7]{HD})}
Let $g:T^{4n-1}\to T^{4n-1}$ be an isometry of $T^{4n-1}$ which is given by formula {\em (\ref{action})} with $A\in SO(4n-2,\Z).$
It extends to $D^2\times T^{4n-2},$ by rotation through the angle $2\pi a$ in the first factor and the extension
has only isolated fixed points.
Suppose that $+1$ is not an eigenvalue of $A.$  
The invariants $\eta_g(0,T^{4n-1})$ are given by
\begin{equation}\label{formula2}
\eta_g(0,T^{4n-1}) = \nu(g)(-1)^{n}\text{cot}(\pi a) \Pi_{i=1}^{2n-1}\text{cot}(\tfrac{\gamma_i}{2})
\end{equation}
where $\nu(g)$ is the number of fixed points of the extension of $g:D^2\times T^{4n-2}\to D^2\times T^{4n-2}$
and $\gamma_i, 1\leq i\leq 2n-1,$ are the rotation angles of $A\in SO(4n-2,\Z).$ The invariants $\nu(g)$ and $\eta_{g}(0,T^{4n-1})$
are independent of the translation $\bar{a}\in\R^{4n-2}/\Z^{4n-2}$ in formula {\em (\ref{action})}.
\end{prop}

\vskip 3mm
\noindent
However, our main result depends from the following observation: 
\newtheorem{cor}{Corollary}
\begin{cor}
Propositions~\ref{prop1} and ~\ref{prop2} are true for $A\in SL(4n-2,\Z).$
\end{cor}
{\bf Proof:} It is well known that any finite order, invertible, integral matrix is conjugate to an orthogonal matrix.
From the third Bieberbach theorem 
(see \cite[Th. 4.1, Chapter I]{charlap}, \cite[Th. 2.1 (3)]{S}, \cite[Th. 3.2.2]{wolf})
we know that an abstract isomorphism between Bieberbach groups can be realized by conjugation within $Aff(n) = GL(n,\R)\ltimes\R^n.$  
Equivalently, it means that flat manifolds with isomorphic
fundamental groups are affine diffeomorphic.
Then it is enough to use theorem 2.4 of \cite{APS} which says, that
the signature $\eta$ invariant is a diffeomorphism invariant.
\vskip 2mm
\hskip 120mm
$\Box$
\section{Main result}
In this section, with the help of formulas (\ref{generalI}) and (\ref{formula2}), we prove that for all seven dimensional
flat manifolds with cyclic holonomy groups, which satisfy condition (\ref{action}) 
the $\eta$ invariant is an integer number. We start with the 3-dimensional case. 
\newtheorem{ex}{Example}
\begin{ex}\label{dimension3}
{\em (See also \cite{LR}.)
Let $M^3$ be a $3$-dimensional, oriented, flat manifold.
There are only six such manifolds. 
The torus and manifolds $M_2, M_3, M_4$, $M_5, M_6,$
with holonomy groups $\Z_2,\Z_3,\Z_4,\Z_6,\Z_2\times\Z_2$ correspondingly,
see \cite[Ch. III]{S} or \cite[Th. 3.5.5]{wolf}.
For holonomy groups $\Z_2$ and $\Z_2\times\Z_2$ the above matrix $A$, has eigenvalues $\pm 1$. 
Hence the $\eta$ invariant is equal to zero. 
In other words $\eta(M_2) = \eta(M_6) = 0.$
Let us consider the case of the
holonomy group $\Z_3.$ Here $M_3 = \R^{3}/\Gamma,$ where 
$$\Gamma = \text{gen}\{g=(B',(1/3,0,0)), (I,(0,1,0)), (I,(0,0,1))\}\subset SL(3,\R)\ltimes\R^3,$$
where
\begin{equation}
B' = \left[
\begin{smallmatrix}
1&0\\
0&B
\end{smallmatrix}\right]\text{with\hskip 2mm}
B = \left[
\begin{smallmatrix}
0&-1\\
1&-1
\end{smallmatrix}\right].
\end{equation}
The action of an element $(A,a)\in SL(3,\R)\ltimes\R^3$ 
is standard: for $x\in\R^3\hskip 3mm (A,a)x = Ax+a.$
We shall use formulas (\ref{generalI}) and (\ref{formula2}). It means that
$\eta(M_3) = -\eta_{g}(0,\R^3/\Z^3)-\eta_{g^2}(0,\R^3/\Z^3),$
where $\Z^3\subset\Gamma$ is the subgroup of translations.
It is easy to see, that the number of fixed points $\nu(g)$
in the formula (\ref{formula2}) is equal to $3.$
In fact, these points are solutions of the matrix equation $BX = X,$ where $X\in \R^{2}/\Z^2.$
It is easy to see that a solution is represented by a three elements, i.e.
$(1/3,2/3), (2/3,1/3), (0,0).$  
Moreover, the rotation angle of the matrix $B$ is equal to $\tfrac{2\pi}{3}$
and $cot(\pi/3) = 1/\sqrt{3}.$
Hence and from formula (\ref{formula2})
$$\eta_{g}(0,\R^3/\Z^3) = -\nu(g)cot(\pi/3)cot(\pi/3) = -1,$$
and
$$\eta_{g^2}(0,\R^3/\Z^3) = -\nu(g^2)cot(2\pi/3)cot(2\pi/3) = -1.$$
Finally
\begin{equation}\label{zet3}
\eta(M_3) = \tfrac{1}{3}(-1 + (-1)) = -\tfrac{2}{3}.
\end{equation}
A similar calculation for the $3$-dimensional flat oriented manifolds with holonomy $\Z_6$ and $\Z_4$ gives
the following version of formula (\ref{zet3}).
For $\Z_6$
$$\eta(M_5) = -\tfrac{1}{6}(2cot^{2}(\pi/6)+6cot^{2}(\pi/3)) = -\tfrac{1}{6}(6+2) = -\tfrac{4}{3}$$
and for $\Z_4$
$$\eta(M_4) = -\tfrac{2}{4}(cot^{2}(\pi/4) + \cot^{2}(3\pi/4)) = -1.$$}
\end{ex}
\newtheorem{theo}{Theorem}
\begin{theo}\label{maintheorem}
Let $T^6$ be any flat six-dimensional torus, and let a cyclic group $G$ act freely on the
Riemannian product $S^1\times T^6$ such that $(x,\bar{x})\to (x+a,A\bar{x}+\bar{a}),$
where $A\in SL(6,\Z)$ descends to $T^6$ and $(a,\bar{a}), (x,\bar{x})\in S^1\times T^6.$
(The action of $G$ satisfies condition (\ref{action}).)  
Let $M^7:= (T^6\times S^1)/G$ be endowed with the induced flat metric. Then $\eta(M^7)\in \Z.$
\end{theo}
{\bf Proof:} 
Let us first assume that the first Betti number $b_1(M^7)\geq 2.$ 
Then from condition (\ref{action}) it follows that the matrix $A$ has an eigenvalue $1,$
see \cite{HS}. Hence, from
Proposition~\ref{prop1},
$\eta(M^7) = 0.$ Further we shall assume that $b_1(M^7) = 1.$ 
\vskip 1mm
From the crystallographic restriction \cite[Proposition 2.1, p. 543]{HH}, the following numbers are
possible for an order of the holonomy group of a seven dimensional flat manifold: 
$2, 3, 4, 5, 6, 7, 8, 9, 10, 12, 14, 15, 18, 20, 24, 30.$
We shall consider these numbers case by case. In each case we shall use the computer 
program CARAT, (see \cite{OPS}, \cite{RL} and \cite{Lu10}) to determine 
the number of flat manifolds which satisfy our assumptions.
\vskip 1mm
\noindent
For a holonomy group of order two, all eigenvalues of the matrix $A$ are equal $\pm 1.$ 
Hence, from Propositions~\ref{prop1},~\ref{prop2} and~\ref{prop3} $\eta(M^7) = 0$
always. 
\vskip 1mm
\noindent
If a holonomy group is equal to $\Z_3$ and $b_1(M^7) = 1,$ we can assume (see \cite[chapter 13]{serre}) 
that the holonomy representation (see (\ref{holonomyrep})) 
is a direct sum (over $\Q$) of 
the trivial representation and 3-times the two-dimensional representations, which we identify with matrices 
$B$ or $B^2$ from the above Example~\ref{dimension3}.  
Hence, we can apply formula (\ref{formula2}). With similar calculations as in the above Example~\ref{dimension3} we have
$\nu(g) = \nu(g^2) = 3^3 = 27,$ where $g\in \pi_1(M^7)$ and $p(g)\neq 1.$ ($p$ was defined on page 1.)
Finally,
$$\eta(M^7) =  \nu(g)cot^{4}(\pi/3)+  \nu(g^2)cot^{4}(2\pi/3) = \tfrac{27}{3}(\tfrac{9}{81}+\tfrac{9}{81}) = 2.$$
\vskip 1mm
\noindent
Let $M^7$ be a flat manifold of dimension $7$ with holonomy group $\Z_4.$ 
From \cite{Lu10}, up to affine diffeomorphism, there are thirteen such
manifolds. For five of them ,the $\eta$ invariant is equal to zero since the matrix $A$ has
an eigenvalue $\pm 1,$ see Proposition~\ref{prop3}. 
Let us consider the case of the diagonal $6\times 6$ matrix
$A=\left[
\begin{smallmatrix}
B&0&0\\
0&B&0\\
0&0&B
\end{smallmatrix}\right],$ where 
$B =\left[
\begin{smallmatrix}
0&-1\\
1&0
\end{smallmatrix}\right].$
Since the set of fixed points of the action of $B$ on $S^1\times S^1$ is equal to $\{(0,0),(1/2,1/2)\},$ 
then $\nu(A) = 8.$\footnote{Here and in what follows $\nu(A) = \nu(g),$ where $g = (A,a),$ see (\ref{action}).} 
Hence
$$\eta(M^7) = \tfrac{8}{4}(cot^{4}(\pi/4) + cot^{4}(3\pi/4)) = 2+2 = 4.$$ 
That is the only case of a flat manifold with holonomy group $\Z_4$
and non-zero $\eta$ invariant, which we can calculate with our methods. 
In Example~\ref{ex2} we define a Bieberbach group 
with a cyclic holonomy group of order four, which does not 
satisfy condition (\ref{action}).
There are seven such manifolds, see \cite{Lu10}.
\vskip 1mm
\noindent
For holonomy groups $\Z_5$ and $\Z_{10}$ it is enough to observe that any flat manifold $M^7$ of dimension $7,$
with such holonomy groups,
either has the first Betti number greater then $1$ or
all matrices $A^k, k\in \Z,$ have eigenvalues $\pm 1.$
In fact, it follows from the crystallographic restriction (see \cite{HH})
that any faithful integral representation of 
the group $\Z_5$ has dimension greater than $3.$
Hence $\eta(M^7) = 0.$ We should add, that a seven dimensional flat manifold with $\Z_5$ holonomy group
and the first Betti number equal to $1$ does not exist, \cite{Lu10}. 
However, there are three such manifolds with holonomy group 
$\Z_{10}.$
\vskip 1mm
\noindent
From \cite{Lu10} there are sixteen isomorphism classes of Bieberbach groups
with holonomy group $\Z_6,$ which are the fundamental groups
of flat $7$-dimensional manifolds with the first Betti number $1.$
All of them satisfy our condition (\ref{action}).
Since in eleven cases the matrix $A$ has eigenvalues $\pm 1,$ then the $\eta$ invariant
is equal to $0,$ cf. Proposition \ref{prop3}.   
Let us assume that $A=\left[
\begin{smallmatrix}
B_{1}&0&0\\
0&B_{2}&0\\
0&0&B_{3}
\end{smallmatrix}\right],$ where 
$B_{i} = D = \left[
\begin{smallmatrix}
0&-1\\
1&1
\end{smallmatrix}\right],$ for $i=1,2,3.$ 
For $A$ the number of fixed points is equal to $1,$ for $A^2$ it is equal to $3^3 = 27.$
Hence, the final formula is the following:
$$\eta(M^7) = \tfrac{1}{6}cot^{4}(\pi/6) + \tfrac{27}{6}cot^{4}(2\pi/6) + $$
$$+ \tfrac{27}{6}cot^{4}(4\pi/6) + \tfrac{1}{6}cot^{4}(5\pi/6) = 4.$$
\vskip 1mm
\noindent
There are still four manifolds to consider. 
For the
first one the matrix $A$ has on the diagonal the matrices $B_1 = D, 
B_2 = B_3 = \bar{B} = \left[
\begin{smallmatrix}
-1&-1\\
1&0
\end{smallmatrix}\right]$ and the $\eta$ invariant is equal to $4.$
For the second case the matrix $A\in GL(6,\Z)$ is not
the diagonal type but is conjugate in $GL(7,\Q)$
to the above matrix. Here the $\eta$ invariant is also equal to $4.$
\vskip 1mm
\noindent  
The last two cases are the following. The matrix $A$ is either the diagonal type
$B_1 = B_2 = D, B_3 = \bar{B}$ or is an integral matrix which 
is conjugate to $A$ in $GL(n,\Q).$
By similar calculation as above, the $\eta$ invariant is equal 
correspondingly to $2$ and $2.$ 
\vskip 3mm
\noindent
For holonomy group $\Z_7$ the matrix  $A=\left[
\begin{smallmatrix}
0&0&0&0&0&-1\\
1&0&0&0&0&-1\\
0&1&0&0&0&-1\\
0&0&1&0&0&-1\\
0&0&0&1&0&-1\\
0&0&0&0&1&-1\\
\end{smallmatrix}\right].$ Then $\nu(A^k) = 7$ and the eigenvalues are equal to 
$$cos(2k\pi/7) + isin(2k\pi/7), k = 1,2,3,4,5,6.$$ Hence,
\footnote{The computations of the trigonometric sums were done with the aid
of a computer and MATHEMATICA version 7.}
$$\eta(M^7) = 2cot(\pi/7)cot(2\pi/7)cot(3\pi/7)(cot(\pi/7)+cot(2\pi/7)-cot(3\pi/7))=2.$$
For holonomy group $\Z_8$ we have eleven isomorphism classes of Bieberbach
groups with the first Betti number one. Six of them do not satisfy condition
(\ref{action}). For two of them, the $\eta$ invariant is equal to zero, since the matrix
$A$ has eigenvalues $\pm 1.$ The last three manifolds have the $\eta$ invariant equal to $2.$
\vskip 1mm
\noindent
For instance, let us present the calculation for the following matrix:
$$\bar{A}=\left[
\begin{smallmatrix}
0&0&0&-1&0&0\\
1&0&0&0&0&0\\
0&1&0&0&0&0\\
0&0&1&0&0&0\\
0&0&0&0&0&1\\
0&0&0&0&-1&0\\
\end{smallmatrix}\right].$$
We have $\nu((\bar{A})^k) = 4$ for $k=1,3,5,7$ and
$$\eta(M^7) = cot(\pi/8)cot(3\pi/8)(cot(\pi/8)-cot(3\pi/8)) = 2.$$
\vskip 1mm
\noindent
For a cyclic group of order $9$ the matrix $$A=\left[
\begin{smallmatrix}
0&0&0&0&0&-1\\
1&0&0&0&0&0\\
0&1&0&0&0&0\\
0&0&1&0&0&-1\\
0&0&0&1&0&0\\
0&0&0&0&1&0\\
\end{smallmatrix}\right].$$
The characteristic polynomial of $A$ is equal to $x^6+x^3+1.$ Moreover $\nu(A^k) = 3$ for $k= 1,2,4,5,7,8$ and
$\nu(A^k) = 27$ for $k = 3,6.$ Hence
$$\eta(M^7) = -\tfrac{2}{3}cot^{2}(\pi/9)cot(2\pi/9)cot(4\pi/9) + \tfrac{2}{3}cot^{2}(2\pi/9)cot(\pi/9)cot(4\pi/9) +$$
$$+ 6cot^{4}(\pi/3) - \tfrac{2}{3}cot^{2}(4\pi/9)cot(\pi/9)cot(2\pi/9) = \tfrac{2}{3} +$$
$$+ \tfrac{2}{3}cot(2\pi/9)cot(\pi/9)cot(4\pi/9)(cot(2\pi/9)-cot(\pi/9)-cot(4\pi/9)) = \tfrac{2}{3}-\tfrac{2}{3} = 0.$$ 
For holonomy group $\Z_{12}$ we have twenty nine manifolds, see \cite{Lu10}.
This is the most non-standard case. 
In seven cases the $\eta$ invariant is zero, because the matrix $A$ has an eigenvalue $\pm 1.$
Moreover, condition (\ref{action}) is not satisfied in ten cases. 
Let us assume, that the matrix 
$A = \left[ 
\begin{smallmatrix}
C_{1}&0&0\\
0&C_{2}&0\\
0&0&C_{3}
\end{smallmatrix}\right],$ where 
$C_{i} = E = \left[
\begin{smallmatrix}
0&-1\\
1&0
\end{smallmatrix}\right],  
C_i = F = \left[
\begin{smallmatrix}
0&-1\\
1&-1
\end{smallmatrix}\right]$ or
$C_i = G = \left[
\begin{smallmatrix}
0&1\\
-1&1
\end{smallmatrix}\right]$
 where
$i=1,2,3.$ 
For example, we shall consider (compare with the case $n=6$) : $C_1 = F,\hskip 2mm C_2 = C_3 = E,$
and $C_1 = C_2 = F,\hskip 2mm C_3 = E.$
It is easy to see, that it is enough to consider in 
the formula (\ref{generalI}) only $A, \hskip 2mm A^{5},\hskip 2mm A^7$ and $A^{11}.$ 
In fact, in all other cases eigenvalues are equal to $\pm 1$ and we can apply 
Propositions \ref{prop1}, \ref{prop2} and \ref{prop3}.
For instance, we have
$$\eta(M^7) = cot(\pi/12)cot(\pi/3)cot^{2}(\pi/4) +  cot(7\pi/12)cot(\pi/3)cot^{2}(\pi/4) +$$  
$$+ cot(2\pi/3)cot(5\pi/12)cot^{2}(\pi/4) +  cot(11\pi/12)cot(2\pi/3)cot^{2}(\pi/4) = 4,$$
or
$$\eta(M^7) = \tfrac{3}{2}(cot(\pi/12)cot^{2}(\pi/3)cot(\pi/4) - cot(7\pi/12)cot^{2}(\pi/3)cot(\pi/4) +$$
$$+ cot^{2}(2\pi/3)cot(\pi/4)cot(5\pi/12) - cot(11\pi/12)cot^{2}(2\pi/3)cot(\pi/4)) =$$
$$= \tfrac{3}{2}(\tfrac{1}{3}cot(\pi/12) + \tfrac{1}{3}cot(5\pi/12) + \tfrac{1}{3}cot(5\pi/12) + \tfrac{1}{3}cot(\pi/12) = 4.$$
\vskip 1mm
\noindent
From \cite{Lu10}, we know that there are six such flat manifolds. For all, the $\eta$ invariant is equal to $4.$
There is still another possibilty: the matrix 
$A =\left[
\begin{smallmatrix}
D&0\\
0&F
\end{smallmatrix}\right],$
where
$$D=\left[
\begin{smallmatrix}
0&-1&0&1\\
1&0&-1&0\\
0&-1&0&0\\
1&0&0&0
\end{smallmatrix}\right]$$  
is the faithful, irreducible rational representation of the group $\Z_{12},$ 
see \cite[p. 234]{BBNWZ}, \cite{OPS}.
Moreover $F=\left[
\begin{smallmatrix}
0&1\\
-1&1
\end{smallmatrix}\right]^{k}$
or
$F=\left[
\begin{smallmatrix}
0&-1\\
1&0
\end{smallmatrix}\right], k = 1,2.$ 
\vskip 1mm
\noindent
From the classification (see \cite{Lu10}) there are six such flat manifolds and 
the $\eta$ invariants are equal 0 in 4 cases and are equal 2 in 2 cases. 
Here, we should mention that one manifold with the $\eta$ invariant
equal to zero is also considered at the end of the paper, see formula (\ref{last}).
\vskip 1mm
\noindent
For holonomy group $\Z_{14},$ $$A=\left[
\begin{smallmatrix}
0&0&0&0&0&1\\
-1&0&0&0&0&1\\
0&-1&0&0&0&1\\
0&0&-1&0&0&1\\
0&0&0&-1&0&1\\
0&0&0&0&-1&1
\end{smallmatrix}\right].$$
The characteristic polynomial of the matrix $A$ is equal to $\tfrac{x^7 + 1}{x + 1}$.
Moreover 
$$\nu(A^k) = \left\{ \begin{array}{ll}
1& \mbox{for $k = 1,3,5,9,11,13$}\\
7& \mbox{for $k = 2,4,6,8,10,12$}\end{array}\right.$$
Finally,
$$\eta(M^7) = cot(\pi/14)cot(3\pi/14)cot(5\pi/14)(\tfrac{1}{7}cot(\pi/14)-\tfrac{1}{7}cot(3\pi/14)-\tfrac{1}{7}cot(5\pi/14)) +$$
$$+ cot(\pi/7)cot(3\pi/7)cot(2\pi/7)(-cot(\pi/7)-cot(2\pi/7)+cot(3\pi/7)) = 0.$$
There is only one flat manifold of this kind, see \cite{Lu10}.
\vskip 2mm
\noindent
For cyclic group of order $15,$ $$A=\left[
\begin{smallmatrix}
0&-1&0&0&0&0\\
1&-1&0&0&0&0\\
0&0&0&0&0&-1\\
0&0&1&0&0&-1\\
0&0&0&1&0&-1\\
0&0&0&0&1&-1
\end{smallmatrix}\right]$$
and $\nu(A^k) = 15$ for $k = 1,2,4,7,8,11,13,14.$ Moreover
$\eta_{g^{k}}(0,T^7)=0,$ for $k=3,5,6,9,10,12,$ where $g\in SO(7)\ltimes \R^7$ denotes an isometry
defined by the formula $$g((x_1,x_2,...,x_7))=(A(x_1,x_2,...,x_6),x_7) + (0,0,...,0,1/15).$$
Summing up, we get
$$\eta(M^7) = 2cot(\pi/3)cot(\pi/5)cot(2\pi/5)(cot(\pi/15) + cot(2\pi/15) +$$
$$+ cot(4\pi/15) - cot(7\pi/15)) = 4.$$
For holonomy group $\Z_{18},$ $$A=\left[
\begin{smallmatrix}
0&0&0&0&0&1\\
-1&0&0&0&0&0\\
0&-1&0&0&0&0\\
0&0&-1&0&0&1\\
0&0&0&-1&0&0\\
0&0&0&0&-1&0
\end{smallmatrix}\right].$$
We have $$\nu(A^k) = \left\{ \begin{array}{lll} 
1& \mbox{for $k = 1,3,5,7,10,11,13,15,17$}\\
3& \mbox{for $k = 2,4,8,10,14,16$}\\
27& \mbox{for $k = 6,12$}\end{array}\right.$$
\vskip 15mm
Hence, 
$$\eta(M^7) = \tfrac{1}{9}cot(\pi/18)cot(7\pi/18)cot(5\pi/18)(-cot(\pi/18) + cot(5\pi/18) - cot(7\pi/18)) +$$
$$+ \tfrac{1}{3}cot(\pi/9)cot(2\pi/9)cot(4\pi/9)(-cot(\pi/9) + cot(2\pi/9) - cot(4\pi/9)) +$$
$$+ \tfrac{1}{9}cot^{4}(\pi/6) + 3cot^{4}(\pi/3) = -1-\tfrac{1}{3}+1+\tfrac{1}{3} = 0.$$
In the cases above, for a holonomy group of order 2, 3, 7, 9, 14, 15 and 18 there exists only one flat manifold
with the first Betti number one.
Let us consider the last three cases of cyclic groups of order $20,24$ and $30.$
We start with the matrix 
$$A=\left[
\begin{smallmatrix}
0&0&0&-1&0&0\\
1&0&0&-1&0&0\\
0&1&0&-1&0&0\\
0&0&1&-1&0&0\\
0&0&0&0&0&-1\\
0&0&0&0&1&0
\end{smallmatrix}\right]$$ of order 20.
From previous facts we have only to consider the following elements:
$A^k, k \in\{1,3,7,9,11,13,17,19\}=S.$ For all other the $\eta$ invariant is equal to zero.
Finally, we have
\vskip 5mm
$$\eta(M^7)=\Sigma_{k\in S}\tfrac{10}{20}(cot(k\pi/20)cot(k\pi/5)cot(2k\pi/5)cot(k\pi/4))=$$
$$=cot(\pi/5)cot(2\pi/5)(cot(\pi/20)+cot(3\pi/20)+cot(7\pi/20)+cot(9\pi/20)) = 4.$$
From \cite{Lu10}, there is another 
flat manifold with holonomy group $\Z_{20},$ and in this case the $\eta$ invariant
is also $4.$ 
The cyclic group of order $24$ is isomorphic to $\Z_8\times\Z_3$ or to $\Z_8\times\Z_6$
and there are only two such flat manifolds, see \cite{Lu10}.
For the first group the matrix 
$$A=\left[
\begin{smallmatrix}
0&0&0&-1&0&0\\
1&0&0&0&0&0\\
0&1&0&0&0&0\\
0&0&1&0&0&0\\
0&0&0&0&0&-1\\
0&0&0&0&1&-1
\end{smallmatrix}\right]$$ and has order $24.$
As in the last case we have only to consider the elements $A^k$ for
$k\in\{1,2,5,7,10,11,13,14,17,19,22,23\}=T.$ Moreover
$$\nu(A^k) = \left\{ \begin{array}{ll}
6& \mbox{for $k = 1,5,7,11,13,17,19,23$}\\
3& \mbox{for $k = 2,10,14,22$}\end{array}\right.$$
\vskip 5mm
We have
$$\eta(M^7)=\Sigma_{k\in \{1,5,7,11,13,17,19,23\}}\tfrac{1}{4}cot(k\pi/24)cot(k\pi/8)cot(3k\pi/8)cot(k\pi/3)+$$
$$\Sigma_{k\in \{2,10,14,22\}}\tfrac{1}{2}cot(k\pi/24)cot(k\pi/8)cot(3k\pi/8)cot(k\pi/3)= 4.$$
For the second group the matrix $A$ is almost the same as the matrix above. 
We only put in the right-down corner, 
the $(2\times 2)$ integral matrix of order $6$ in place of the matrix of order $3.$ 
The $\eta$ invariant is equal to $0.$ 
In the last case of cyclic group of order 30 there are three manifolds with
the following matrices
$$A_1=\left[
\begin{smallmatrix}
0&0&0&-1&0&0\\
1&0&0&-1&0&0\\
0&1&0&-1&0&0\\
0&0&1&-1&0&0\\
0&0&0&0&0&1\\
0&0&0&0&-1&1
\end{smallmatrix}\right],
A_2 = \left[
\begin{smallmatrix}
0&0&0&1&0&0\\
-1&0&0&1&0&0\\
0&-1&0&1&0&0\\
0&0&-1&1&0&0\\
0&0&0&0&-1&-1\\
0&0&0&0&1&0
\end{smallmatrix}\right]$$ and
$$A_3 = \left[
\begin{smallmatrix}
0&0&0&1&0&0\\
-1&0&0&1&0&0\\
0&-1&0&1&0&0\\
0&0&-1&1&0&0\\
0&0&0&0&0&1\\
0&0&0&0&-1&1
\end{smallmatrix}\right].$$ 
We shall present calculatins for $A_1$, the other cases are similar.
As in the cases above, we have only to consider elements 
$$A_{1}^{k}, k\in\{1,2,4,7,8,11,12,13,14,16,17,19,22,23,26,28,29\}=R.$$
We have
$$\nu(A^k) = \left\{ \begin{array}{ll}
5& \mbox{for $k = 1,7,11,13,17,19,29$}\\
15& \mbox{for $k = 2,4,8,14,16,22,26,28$}\end{array}\right.$$
\vskip 5mm
We have
$$\eta(M^7)=
\Sigma_{k\in\{1,7,11,13,17,19,29\}}\tfrac{5}{30}cot(k\pi/30)cot(2k\pi/10)cot(4k\pi/10)cot(k\pi/6)+$$
$$+\Sigma_{k\in\{2,4,8,14,16,22,26,28\}}\tfrac{15}{30}cot(k\pi/30)cot(2k\pi/10)cot(4k\pi/10)cot(k\pi/6)= 0.$$
For the other two manifolds with
holonomy group of order $30,$ the $\eta$ invariant is equal to $4$ and $0.$
\vskip 2mm
\hskip 120mm
$\Box$
\vskip 2mm
\noindent
Let us present a final table.\begin{center}
\begin{tabular}{|c|c|c|c|c|l|}
\hline
holonomy & number of & number & \multicolumn{2}{|c|} {number of cases not calculated} & values of $\eta$ \\
\cline{4-5}
&  manifolds & with $b_1 = 1$ & all & $b_1 = 1$ &  \\
\hline
$\Z_2$ & 15 & 1 & 0 & 0 & 0\\
\hline
$\Z_3$ & 6 & 1 & 0 & 0& 0, 2\\
\hline
$\Z_4$ & 87 & 13 & 37 & 7 & 0, 4\\
\hline
$\Z_5$ & 2 & 0 & 0 & 0 & 0\\
\hline
$\Z_6$ & 74 & 16 & 0 & 0 & 0, 2, 4\\
\hline
$\Z_7$ & 1 & 1 & 0 & 0 & 0\\
\hline
$\Z_8$ & 24 & 11 & 10 & 6 & 0, 2\\
\hline
$\Z_{9}$ & 1 & 1 & 0 & 0 & 0\\
\hline
$\Z_{10}$ & 12 & 3 & 0 & 0 & 0\\
\hline
$\Z_{12}$ & 89 & 29 & 27 & 10 & 0, 2, 4\\
\hline
$\Z_{14}$ & 1 & 1 & 0 & 0& 0\\
\hline
$\Z_{15}$ & 1 & 1 & 0 & 0 & 4\\
\hline
$\Z_{18}$ & 1 & 1 & 0 & 0 & 0\\
\hline
$\Z_{20}$ & 2 & 2 & 0 & 0 & 4\\
\hline
$\Z_{24}$ & 2 & 2 & 0 & 0 & 0, 4\\
\hline
$\Z_{30}$ & 3 & 3 & 0 & 0 & 0, 4\\
\hline
\end{tabular}
\end{center}
\vskip 10mm
\begin{ex}\label{ex2}
{\em Let $\R^7/\Gamma$ be the seven dimensional flat oriented manifold with the fundamental group $\Gamma\subset SO(7)\ltimes\R^7$
generated by
$$\{(A,(1/2,0,...,0)),(I,e_i)\},$$
where $e_i, i=1,2,...,7$ is a standard basis of $\R^7$ and 
$$A=\left[
\begin{smallmatrix}
0&1&0&0&0&0&0\\
1&0&0&0&0&0&0\\
0&0&0&1&0&0&0\\
0&0&1&0&0&0&0\\
0&0&0&0&0&-1&0\\
0&0&0&0&1&0&0\\
0&0&0&0&0&0&1
\end{smallmatrix}\right].$$ It is easy to see that the manifold $\R^7/\Gamma$
has holonomy $\Z_4$ and does not satisfy the condition (\ref{action}).}
\end{ex}
\section{Holonomy group as a subgroup of $SO(n,\Z)$}
This part is a modified and refreshed version of some results of \cite{HD}.
We shall present the formula (\ref{formula2}) 
under the assumption that the image of the holonomy representation
(see (\ref{holonomyrep})) is a subgroup of $SO(n,\Z).$
In this case (see \cite{HD}) a method is given for finding the number $\nu(g)$ of
fixed points of $g,$ where $g$ is an isometry of the torus. 
\vskip 1mm
\noindent
Let $e_i, 0\leq i\leq 4n-2$ be a standard basis of the space $\R^{4n-2}.$
Since $A\in SO(4n-2,\Z)$ one has $A(e_i)=\pm e_{j(A,i)}, i\neq j(A,i)$ for each $i.$ 
We denote by $\sigma(A)$ the element of $SO(4n-2,\Z)$ defined
by $\sigma(A)(e_i) = e_{j(A,i)}.$ 
Then $\sigma(A)$ is a permutation matrix and we may decompose $\sigma(A) = \sigma_1(A)\sigma_2(A)...\sigma_l(A)$
into disjoint cycles.
Summing up we have:
\begin{prop}\label{prop3} {\em (\cite[Proposition 4.9]{HD})}
Let $g:T^{4n-1}\to T^{4n-1}$ be an isometry given by the formula (\ref{action}) with $a\neq 1.$ 
Then $\eta_g(0,T^{4n-1}) = 0$ if $A$ has an eigenvalue equal to $+1$ or $-1.$
Otherwise
\begin{equation}\label{formula3}
\eta_g(0,T^{4n-1}) = 2^l(-1)^{n}\text{cot}(\pi a)\Pi_{i=1}^{2n-1}\text{cot}(\frac{\gamma_i}{2})
\end{equation}
where
\vskip 1mm
(i) The angles $\gamma_i$ are the rotation angles of the orthogonal matrix $A.$
\vskip 1mm
(ii) The integer $l$ is the number of distinct cycles in the decomposition of the permutation matrix
$$\sigma(A) = \sigma_1(A)\sigma_2(A)....\sigma_l(A),$$ corresponding to $A.$
\end{prop}

\vskip 2mm
\hskip 120mm
$\Box$
\vskip 2mm
As an immediate corollary we have a special version of the formula (\ref{general}).
\begin{prop}\label{prop4}{\em (\cite[Proposition 4.12]{HD})}
We keep the above notations. The corresponding eta invariant is given by
\begin{equation}\label{finalformula}
\eta(0,Y) = 
\tfrac{1}{\mid G\mid}\Sigma'(2^l(-1)^{n}\text{cot}(\pi a)\Pi_{i=1}^{2n-1}\text{cot}(\tfrac{\gamma_i}{2}))
\end{equation}
where the symbol $\Sigma'$ means summation over the group elements $g$ whose associated $A$ has no eigenvalue equal to $\pm 1.$ 
Furthermore
$\gamma_i$ are the rotations angles of $A\in SO(4n-2,\Z),$ and  $l$ is the number of distinct cycles of $A$ as above.
\end{prop}
\vskip 1mm 
Let us present two families of flat manifolds which satisfy the above assumptions.
\begin{ex} 
{\em 1. Let $M^n$ be an oriented flat manifold with holonomy group $(\Z_2)^{n-1},$
so called Hantzsche-Wendt manifold. From \cite{HW} it follows that
$n$ is an odd number and the image of the holonomy representation of $\pi_1(M^n)$
(see (\ref{holonomyrep})) in $GL(n,\Z)$ consists of diagonal matrices with $\pm 1$
on the diagonal. Hence, from Proposition~\ref{prop3}, $\eta(M^n) = 0.$
We should add, that the same is true for any oriented flat manifold of dimension $n$
with holonomy group $(\Z_2)^k, 1\leq k\leq n-2$ whose fundamental group has the
holonomy representation with a diagonal image in the above sense.}
\end{ex}
\vskip 2mm
\begin{ex}
{\em Let us recall the fundamental group $\Gamma$ of a flat manifold of dimension $4n-1$ with holonomy group $\Z_{2(4n-2)}.$
From Bieberbach theorems (see \cite{S}) $\Gamma$ defines the short exact sequence of groups.
\begin{equation}\label{cyslicb}
0\rightarrow \Z^{4n-1}\rightarrow\Gamma\rightarrow\Z_{2(4n-2)}\rightarrow 0.
\end{equation}
Moreover $\Gamma\subset SO(4n-1)\ltimes\R^{4n-1},$ is generated by 
the translation subgroup $\Gamma\cap\{I\}\times\R^{4n-1}\simeq\Z^{4n-1}$ and
an element $(A',(0,...,0,\tfrac{1}{2(4n-2)})).$ Here 
\begin{equation}
A' = \left[
\begin{smallmatrix}
A&0\\
0&1
\end{smallmatrix}\right]
\end{equation}
and 
\begin{equation}
A = \left[
\begin{smallmatrix}
0&0&...&0&-1\\ 
1&0&...&0&0\\
0&1&0&...&0\\
...........\\
0&0&...&1&0
\end{smallmatrix}\right]
\end{equation}
is $(4n-2)\times(4n-2)$ integral orthogonal matrix of order $2(4n-2).$
Using the formula (\ref{finalformula}) 
from Proposition \ref{prop4} for a flat manifold $\R^{4n-1}/\Gamma = M^{4n-1}$ we have
\begin{equation}\label{mainformula}
\eta(M^{4n-1}) = \tfrac{(-1)^n}{8n-4}\Sigma_{k=1}^{8n-5}2^{l_{k}}(\text{cot}(\frac{k\pi}{2(4n-2)})\Pi_{l=1}^{2n-1}\text{cot}(\frac{k(2l-1)\pi}{2(4n-2)})),
\end{equation}
where $l_k$ is the number of distinct cycles of $A^k,$ see Proposition~\ref{prop3}.
The characteristic polynomial of the matrix $A$ is 
equal to det$(A-(\lambda)I) = \lambda^{4n-2}+1.$ Hence,
it is easy to calculate that $A^{2n-1} = -I,$ and $\pm 1$ is 
a root of the characteristic polynomial of $A^2.$ 
It is to see that, for $1\leq m < n$ and $k=2r+1$ 
\begin{equation}\label{elem1}
cot(\tfrac{k(m-1)\pi}{2(4n-2)}) = 
cot(\tfrac{k\pi}{2} - \tfrac{k(2m-1)\pi}{2(4n-2)}) = tg(\tfrac{k(2m-1)\pi}{2(4n-2)})
\end{equation}
and
\begin{equation}\label{elem2}
cot(\tfrac{k(2n-1)\pi}{2(4n-2)}) = cot(\tfrac{k\pi}{4}) = (-1)^r.
\end{equation}
\vskip 1mm
\noindent
Using elementary formulas
$$cot(\tfrac{\pi}{2}-\alpha) = tg(\alpha),\hskip 2mm cot(\pi-\alpha) = - ctg(\alpha)$$
and the above equations (\ref{elem1}), (\ref{elem2}) we obtain
\begin{equation}\label{wzor}
\eta(M^{4n-1}) = \tfrac{(-1)^n}{8n-4}\Sigma_{r=0}^{4n-3}2^{l_{r}}(-1)^{r}
\text{cot}(\tfrac{(2r+1)\pi}{2(4n-2)}).
\end{equation}
\vskip 1mm
\noindent
Here $l_r$ is the number of distinct cycles of $A^{(2r+1)}.$
\vskip 2mm
\noindent
For $n=2,$ the components of the above formula are non zero only for $k = 1,3,5,7,9,11.$
It is easy to see that $l_{0} = l_{2} = l_{3} = l_{5}$
and $l_{1} = l_{4} = 3.$ Summing up
\begin{equation}\label{last}
\eta(M^7) =  \tfrac{2}{12}cot(\tfrac{\pi}{12}) - \tfrac{8}{12} + \tfrac{2}{12}tg(\tfrac{\pi}{12}) + 
\tfrac{2}{12}tg(\tfrac{\pi}{12}) - \tfrac{8}{12} +
\end{equation} 
$$+ \tfrac{2}{12}cot(\tfrac{\pi}{12}) = \tfrac{1}{3}(cot(\tfrac{\pi}{12}) + tg(\tfrac{\pi}{12})) - 
\tfrac{4}{3} = \tfrac{4}{3} - \tfrac{4}{3} = 0\in\Z.$$}
\end{ex}

\end{document}